\documentclass[titlepage,11pt,fleqn]{article}
\usepackage{amsthm}
\newtheorem{thm}{}[section]

\usepackage{latexsym}
\usepackage{amsfonts}

\usepackage[leqno]{amsmath}

\makeatletter
\newcommand{\leqnomode}{\tagsleft@true}
\newcommand{\reqnomode}{\tagsleft@false}
\makeatother

\usepackage{enumerate}
\usepackage{tikz}
\usetikzlibrary{graphs}
\usepackage{float}
\usepackage{amssymb}

\newcommand\blackslug{\hbox{\hskip 1pt \vrule width 4pt height 8pt depth 1.5pt
        \hskip 1pt}}
\newcommand\bbox{\hfill \quad \blackslug \bigbreak}
\def\dd{\hbox{-}}
\def\cc{\hbox{-}\cdots\hbox{-}}
\def\ll{,\ldots,}
\newcommand{\sset}[1]{\left\{#1\right\}}
\newcommand{\Proof}{\noindent{\bf Proof.}\ \ }

\title{Finding an induced path that is not a shortest path}
\author{
Eli Berger\thanks{Supported by Israel Science Foundation Grant 100004639
and Binational Science Foundation USA-Israel Grant 100005728.}\\
University of Haifa\\ \\
Paul Seymour\thanks{Supported by AFOSR grant A9550-19-1-0187 and NSF
grant DMS-1800053.}\\
Princeton University, Princeton, NJ 08544\\ \\
Sophie Spirkl\thanks{Current address: Princeton University, Princeton, NJ08544. This material is based upon work supported by the National Science Foundation under Award No. DMS-1802201.}\\
Rutgers University, Piscataway, NJ 08854}

\date{July 18, 2019; revised \today}

\begin{document}
\maketitle
\begin{abstract}
We give a polynomial-time algorithm that, with input a graph $G$ and two vertices $u,v$ of $G$, decides
whether there is an induced $uv$-path that is longer than the shortest $uv$-path.
\end{abstract}

\section{Introduction}

All graphs in this paper are finite and simple. For a graph $G$ and $u, v\in V(G)$, the \emph{$G$-distance} $d_G(u,v)$ ($d(u,v)$ when 
there is no danger of confusion) is the number of edges in a shortest $uv$-path in $G$; let $d(u, v) = \infty$ if there is no 
such path. Let $P$ be an induced $uv$-path. 
The \emph{length} of $P$ is the number of edges of $P$. We call $P$ a \emph{non-shortest $uv$-path ($uv$-NSP)} if the length of $P$ 
is more than $d(u,v)$. 

Given a graph $G$ and $u, v \in V(G)$ we consider the question of whether there are two induced $uv$-paths of different lengths, or 
equivalently, whether there is a $uv$-NSP. Deciding this in polynomial time is surprisingly non-trivial. (It is important that we
want induced paths; if we just want paths of different lengths, the question is much easier.)
Our main result is the following:

\begin{thm} \label{thm1} \label{thm:long}
There is an algorithm that, given a graph $G$ and $u, v \in V(G)$, decides whether there is a $uv$-NSP in time $O(|G|^{16})$.
\end{thm}

A step in the proof has the following consequence which may also be of interest:

\begin{thm} \label{thm:longbutshort}
For fixed $k$, there is a polynomial-time algorithm that, given a graph $G$ and $u, v \in V(G)$, decides whether there is an 
induced path between $u$ and $v$ 
in $G$ of length exactly $d(u,v) + k$. 
\end{thm}

We prove \ref{thm:longbutshort} in section 2, and \ref{thm:long} in section \ref{sec:nsp}.
Many variants of finding pairs of induced paths have been considered previously; for instance

\begin{thm}[Bienstock \cite{bienstock}]
The following problems are $NP$-hard: 
\begin{itemize}
\item Given $u, v \in V(G)$, decide whether there is an induced $uv$-path of odd (even) length. 
\item Given $u, v \in V(G)$, decide whether there are two induced $uv$-paths $P_1$ and $P_2$ with no edges between $V(P_1) \setminus \{u,v\}$ and $V(P_2) \setminus \{u,v\}$.
\end{itemize}
\end{thm}



Here are two more NP-hardness results, that are new as far as we know, but for reasons of space we omit the proofs:
\begin{thm}
\label{twopaths}
The following problem is NP-hard: 
\begin{itemize}
    \item Input: A graph $G$ and $u,v \in V(G)$.
    \item Output: ``Yes'' if there exist two induced $uv$-paths $P$ and $Q$ such that there are no edges between 
$V(P) \setminus \sset{u,v}$ and $V(Q) \setminus \sset{u,v}$, and $P$ is a shortest $uv$-path; and ``No'' otherwise. 
\end{itemize}
\end{thm}

This is in contrast with \ref{lem:manyshortpaths}, which implies that the problem is polynomial-time solvable if both $P$ and $Q$ are both 
required to be shortest paths (or at most a fixed constant amount longer than a shortest path). 
In view of \ref{thm1}, it is natural to ask: 

\begin{thm}
For fixed $k > 1$, is there a polynomial-time algorithm that, given a graph $G$ and $u, v \in V(G)$, decides 
whether there is an induced $uv$-path $P$ in $G$ of length at least $d(u,v) + k$?
\end{thm}

This remains open, even for $k=3$ (the algorithm of this paper does the case $k=1$, and can be adjusted to do the case $k=2$).
It is necessary to fix $k$, because of the following: 

\begin{thm}
\label{twiceaslong}
The following problem is NP-hard: 
\begin{itemize}
    \item Input: A graph $G$ and $u,v \in V(G)$.
    \item Output: ``Yes'' if there exists a $uv$-NSP of length at least $2d_G(u,v)$ and ``No'' if there is no such path. 
\end{itemize}
\end{thm}

\section{Dynamic programming}

A {\em path forest} means a graph in which every component is a path (possibly of length zero); and a {\em path forest in $G$}
means an induced subgraph of $G$ that is a path forest. (Thus it consists of a set of induced paths of $G$, pairwise vertex-disjoint and
with no edges of $G$ joining them.)

Let $V_1\ll V_n$ be pairwise disjoint subsets of $V(G)$, with union $V(G)$, such that for all $i,j\in \{1\ll n\}$, if $j\ge i+2$
then there are no edges between $V_i$ and $V_j$. We call $(V_1\ll V_n)$ an {\em altitude}. We are given a graph $G$ and an altitude 
$(V_1\ll V_n)$ in $G$, and we need to test whether there is a path forest in $G$ with certain properties, that contains only a 
bounded number of vertices from each $V_i$. We shall see that this can easily be solved with dynamic programming.

Let $X\subseteq V(G)$, and let $H, H'$ be path forests in $G$. We say they are {\em $X$-equivalent} if
\begin{itemize}
\item $V(H)\cap X=V(H')\cap X$;
\item $H$, $H'$ have the same number of components; and
\item for each component $P$ of $H$, there is a component $P'$ of $H'$ with the same ends and same length as $P$.
\end{itemize}
This is an equivalence relation.

Again, let $X\subseteq V(G)$. A path forest $H$ is {\em $h$-restricted} in $G$ relative to $X$ if $|V(H)\cap X|\le h$,
and there are at most $h$ components of $H$ that have no end in $X$.
Now let $(V_1\ll V_n)$ be an altitude in $G$. A path forest $H$ is {\em $h$-narrow} (with respect to $(V_1\ll V_n)$) if
for $1\le i\le n$, $H[V_i\cup\cdots\cup V_n]$ is $h$-restricted in $G[V_i\cup\cdots\cup V_n]$ with respect to $V_i$.

Let $1\le i\le n$. Let $\mathcal{C}_i$ be the set of all equivalence classes, under $V_i$-equivalence, that contain
a path forest in $G[V_i\cup\cdots\cup V_n]$ that is $h$-narrow with respect to $(V_i\ll V_n)$. Algorithmically, 
we may describe $\mathcal{C}_i$ by explicitly storing such a path forest. 

We observe:
\begin{thm}\label{update} If $h$ is fixed, with $G$, $V_1\ll V_n$ as above, for $1\le i<n$ we can compute
$\mathcal{C}_i$ from a knowledge of $\mathcal{C}_{i+1}$ in polynomial time.
\end{thm}
\Proof There are only polynomially many equivalence classes in $\mathcal{C}_{i+1}$. (This is where we use the condition that at most 
$h$ components of $H$ have no end in $X$, in the definition of ``$h$-restricted''.) For each one, take a representive member $H'$ say.
There are only polynomially many induced subgraphs $J$ of the graph $G[V_i\cup V_{i+1}]$ such that $V(J)\cap V_{i+1}=V(H')\cap V_{i+1}$
and $|V(J)\cap V_i|\le h$. For each such $J$, check whether $H'\cup J$ is $h$-narrow in $G[V_i\cup\cdots\cup V_n]$ 
with respect to $(V_1\ll V_n)$, and if so record 
its equivalence class under $V_i$-equivalence. To see that every member of $\mathcal{C}_i$ is recorded, observe that 
if $H$ is a path forest in $G[V_i\cup\cdots\cup V_n]$ that is $h$-narrow with respect to $(V_i\ll V_n)$, then
$H\setminus V_i$ is a path forest in $G[V_{i+1}\cup\cdots\cup V_n]$ that is $h$-narrow with respect to $(V_{i+1}\ll V_n)$; and
if $H'$ is another member of the equivalence class in $\mathcal{C}_{i+1}$ that contains $H\setminus V_i$, then its union with
$J=H[V_i\cup V_{i+1}]$ is $h$-narrow with respect to $(V_{1}\ll V_n)$ and $V_i$-equivalent to $H$.  
This proves \ref{update}.~\bbox

We deduce:
\begin{thm}\label{genalg}
For all fixed $h\ge k\ge 0$, there is a polynomial-time algorithm that, given pairs $(s_1,t_1)\ll (s_r,t_r)$ of a graph $G$, 
and integers $n_1\ll n_r\ge 0$, 
and an altitude $(V_1\ll V_n)$ in $G$, computes whether there is a path forest in $G$, $h$-restricted with respect to $(V_1\ll V_n)$,
with $r$ components, where the $i$th component has ends $s_i, t_i$ and has length $n_i$.
\end{thm}
\Proof 
First compute $\mathcal{C}_n$; then $n-1$ applications of \ref{update} allow us to compute $\mathcal{C}_1$, and from $\mathcal{C}_1$
we can read off the answer.~\bbox

This implies \ref{thm:longbutshort2}, which we restate:
\begin{thm} \label{thm:longbutshort2}
For fixed $k$, there is a polynomial time algorithm that, given a graph $G$ and $u, v \in V(G)$, decides whether there is an induced path 
between $u$ and $v$
in $G$ of length exactly $d(u,v) + k$. 
\end{thm}
We may assume that $G$ is connected. For each $i\ge 0$, let $V_i$ be the set of vertices with distance exactly $i$ from $u$. Then $(V_1\ll V_n)$
is an altitude, where $n$ is the largest $i$ with $V_i\ne \emptyset$.
Let $P$ be an induced $uv$-path of length $d(u,v)+k$. Then, 
for all $i \in \{1, \dots, d(u,v)\}$, $P$ contains a vertex $x$ with $d(x,v) = i$. Consequently, for all $i \in \mathbb{N}_0$, 
$P$ contains at most $k+1$ vertices with distance exactly $i$ from $v$. So $P$ is $(k+1)$-narrow with respect to $(V_1\ll V_n)$,
where $n$ is the largest $i$ with $V_i\ne \emptyset$. Hence \ref{genalg}, with $r=1$ and $n_1=d(u,v) + k$, 
will detect a path in the same $V_1$-equivalence class.~\bbox

Similarly, by trying all possibilities for $n_1\ll n_r$, we obtain  
\begin{thm} \label{lem:manyshortpaths} For fixed $h$ and $r$, there is a polynomial-time algorithm with the following specifications, where
$V_i$ is the set of vertices with distance exactly $i$ from $v$:
\begin{itemize}
    \item Input: A graph $G$, $v \in V(G)$ and $r$ pairs $(s_1, t_1), \dots, (s_r, t_r) \in V(G)$.
    \item Output: A path forest $H$ of $G$ with $r$ components $P_1\ll P_r$, such that for each $i$, $P_i$ has ends $s_i, t_i$
and $|V(H) \cap V_j| \leq h$ for all $j \in \mathbb{N},$ or a determination that no such path forest exists.
\end{itemize}
\end{thm}

\section{Finding an induced non-shortest path} \label{sec:nsp}

In this section, we prove \ref{thm:long}. We start with some definitions. 
A vertex $x \in V(G)$ is \emph{$uv$-straight} if $d(u,x)+d(x,v) = d(u,v)$. Let $G$ be a graph, and $u, v \in V(G)$. Let 
$F$ be the set of $uv$-straight vertices. For $i \in \{0, \dots, d(u,v)\}$, let 
$V_i = \{x \in F: d(u,x) = i\}$; we call $V_i$ the \emph{$uv$-layer of height $i$}, and we say its elements have {\em height} $i$; and we call the sequence 
$V_0, \dots, V_{d(u,v)}$ the \emph{$uv$-layering} of $G$. It follows that for $i, j \in \{0, \dots, d(u,v)\}$ with $|i-j| \geq 2$, 
there are no edges between $V_i$ and $V_j$, and moreover, for $i \in \{1, \dots, d(u,v)-1\}$, every vertex in $V_i$ has a neighbour 
in $V_{i-1}$ and in $V_{i+1}$. 

We call a path $Q$ with $V(Q) \subseteq F$ {\em monotone} (leaving the dependence on $u,v$ to be understood) if 
$|V(Q) \cap V_i| \leq 1$ for all 
$i \in \{0, \dots, d(u,v)\}$ (and therefore $Q$ is induced); and it follows that the vertices of $Q$ are in $|V(Q)|$ $uv$-layers of consecutive heights. 
For every vertex $x \in F$, there is a monotone $xu$-path intersecting precisely $V_0, \dots, V_{d(u,x)}$ and a monotone 
$xv$-path intersecting precisely $V_{d(u,x)}, \dots, V_{d(u,v)}$, and from the definition of $uv$-monotonicity, it follows that 
both of these paths are shortest paths. If $K\subseteq V(G)$, $N(K)$ or $N_G(K)$ denotes the set of all vertices in $V(G)\setminus K$
that have a neighbour in $K$.

Conveniently, in order to solve \ref{thm1} it is enough to handle the case when all vertices are $uv$-straight, because of the next result.

\begin{thm} \label{lem:components}
There is a polynomial-time algorithm with the following specifications: 
\begin{itemize}
    \item Input: A graph $G$ and $u,v \in V(G)$.
    \item Output: Either a $uv$-NSP, or a graph $G'$ with $u,v\in V(G')\subseteq V(G)$ such that $G'$ has a $uv$-NSP if and only if 
$G$ has a $uv$-NSP, and such that every vertex of $G'$ is $uv$-straight in $G'$. 
\end{itemize}
\end{thm}
\Proof
Let $G$ be a graph, and $u, v \in V(G)$. We compute the set $F$ of $uv$-straight vertices, 
and the $uv$-layering $V_0, \dots, V_{d(u,v)}$ of $G$. We may assume that $V(G) \setminus F \neq \emptyset$, for otherwise $G, u, v$ is the desired output. 

Compute the vertex set $K$ of a connected component of $G \setminus F$.  Suppose first that $N(K)$ contains non-adjacent vertices $x,y$
with $d(u, x) < d(u, y)$, and choose $x, y$ such that $d(u,y) - d(u,x)$ is maximum. Let $i=d(u,x)$ and $j=d(u,y)$. 
It follows that no vertex in $V_0, \dots, V_{i-1}$ has a neighbour in $K$ (for otherwise such a vertex contradicts the choice of $x$); 
and similarly, no vertex in $V_{j+1}, \dots, V_{d(u,v)}$ has a neighbour in $K$. Now let $P_1$ be a monotone $xu$-path, let $P_2$ be a 
monotone $yv$-path, and let $Q$ be an induced $xy$-path with interior in $K$. It follows that the concatenation $P_1\dd Q\dd P_2$ is an 
induced $uv$-path; and since $V(Q) \cap K \neq \emptyset$, it follows from the definition of $K$ and $F$ that $P_1\dd Q\dd P_2$ is a 
$uv$-NSP, and we can find such a path in polynomial time. 

Thus we may assume that $N(K)$ is contained in $V_i \cup V_{i+1}$ for some $i \in \sset{0, \dots, d(u,v)-1}$, and $N(K)\cap V_i$ 
is complete to $N(K) \cap V_{i+1}.$ Let $H$ be obtained from $G$ by deleting $K$ and adding edges to make $N(K)$ a clique.  
We claim that $H$ has a $uv$-NSP if and only if $G$ does. 

Suppose first that $P$ is a $uv$-NSP of $G$. Since $N(K)$ is a clique of $H$, there is a $uv$-path of $H$ with vertex set a 
subset of $V(P)$; let $Q$ be the shortest such path.
We claim that $Q$ is a $uv$-NSP of $H$. If $V(P) = V(Q)$, this follows from the choice of $P$. Otherwise, $Q$ contains 
an edge $e$ in $E(H) \setminus E(G)$. Since $e$ connects two vertices at the same distance from $u$, it follows that every
induced $uv$-path containing $e$ is a $uv$-NSP of $H$, as claimed, and so $H$ has a $uv$-NSP.  

Now suppose that $Q$ is a $uv$-NSP of $H$. If $Q$ does not contain an edge in $E(H) \setminus E(G)$, then $Q$ is a $uv$-NSP of $G$, so
we assume that $Q$ contains such an edge. 
Since $N(K)$ is a clique of $H$, it follows that $Q$ contains exactly two vertices $x, y \in N(K)$, and $xy \not\in E(G)$. 
Let $P$ be obtained from $Q$ by replacing $xy$ by an induced $xy$-path with interior in $K$. Then $P$ is a $uv$-NSP of $G$, since $P$ 
contains a vertex of $K$. This proves that $H$ has a $uv$-NSP if and only if $G$ does. 

By repeating this procedure for all components of $G \setminus F$, we either find a $uv$-NSP, or the desired graph $G'$.\bbox

\begin{thm}\label{thm:k0}
There is a polynomial-time algorithm with the following specifications: 
\begin{itemize}
    \item Input: A graph $G$ and $u,v \in V(G)$ such that every vertex of $G$ is $uv$-straight. 
    \item Output: A $uv$-NSP in $G$, or a determination that none exists.
    \item Running time: $O(|G|^{16})$.
\end{itemize}
\end{thm}
\Proof
For $i \in \sset{0, \dots, d(u,v)}$, let $V_i = \sset{x \in V(G): d(x, u) = i}$, and for each vertex $x$, let $h(x)$ be its height. 
Let $P$ be a shortest $uv$-NSP in $G$ (if one exists). We will prove some properties of $P$ that will make it easier to find $P$. 

Let $P_u$ be the longest monotone subpath of $P$ containing $u$, and let $P_v$ be the longest monotone subpath of $P$ containing $v$. 
Let $s$ denote the endpoint of $P_u$ that is not $u$, and let $t$ denote the endpoint of $P_v$ that is not $v$. It follows that $P_u$ and $P_v$ are disjoint, for otherwise $P$ is monotone, contrary to the choice of $P$. 
\\
\\
(1) {\em $V(P) \setminus V(P_v)$ does not contain a vertex $x$ with $h(x) > h(s)$, and $V(P) \setminus V(P_u)$ does not contain a 
vertex $x$ with $h(x) < h(t)$.}
\\
\\
Let $x \in V(P) \setminus V(P_v)$ be chosen with $h(x)$ maximum, breaking ties by choosing the vertex closest to $u$ along $P$. 
Let $Q$ be a monotone $xv$-path, and let $P'$ be the subpath of $P$ from $u$ to $x$. Let $Q'$ denote the concatenation of $P'$ and $Q$. 
We claim that $Q'$ is shorter than $P$. This follows since the subpath of $P$ from $x$ to $v$ is not monotone (because $x \not\in V(P_v)$),
and the subpath of $Q'$ from  $x$ to $v$ is monotone. Since $P$ is a shortest $uv$-NSP, it follows that $Q'$ is not a $uv$-NSP, 
and hence $Q'$ is monotone. In particular, $P'$ is monotone. Thus $V(P') \subseteq V(P_u)$. From the choice of $x$, it follows that 
$P' = P_u$; and so $u=s$. From the choice of $x$, and from the symmetry between $u$ and $v$, this proves (1).

\bigskip

Since $P$ is not monotone, (1) immediately implies that $h(s) \ge h(t)$. 
\\
\\
(2) {\em For fixed $k$, if $h(s)-h(t)\le k$, then we can find a $uv$-NSP in polynomial time (depending on $k$).}
\\
\\
It suffices to prove (2) when $h(s)-h(t) = k$; then we obtain the desired algorithm by applying the statement for $k' = 0, \dots, k$. 

Let $xy \in E(G)$ with $h(y) = h(x)+1$,
and let $v_1 \cc v_{k+2}$ be a $(k+2)$-vertex path with $h(v_i) = h(y) + i - 1$ for $1\le i\le k+2$,
such that $v_1$ is nonadjacent to $x$, and $v_i$ is nonadjacent to $x,y$ for $2\le i\le k+2$.
For all such choices of $x, y, v_1\ll v_{k+2}$, we proceed as follows: 
\begin{itemize}
    \item Let $Q_u$ be a monotone path from $x$ to $u$, and let $Q_v$ be a monotone path from $v_{k+2}$ to $v$. 
    \item We delete all vertices and neighbours of $V(Q_u) \cup V(Q_v) \cup \sset{x} \cup \sset{v_2, \dots, v_{k+2}}$ 
except for $y$ and $v_1$ from $G$. Let $H$ denote the graph we obtain by these deletions. 
    \item We check if $H$ contains an induced path $Q$ from $v_1$ to $y$. If so, we return the concatenated path 
$$Q' = u\dd Q_u\dd x\dd y \dd Q\dd v_1\dd v_2\cc v_{k+2}\dd Q_v\dd v.$$
\end{itemize}

First, we claim that if this returns a path $Q'$, then $Q'$ is a $uv$-NSP. From the construction of $H$, it follows that $Q'$ is an 
induced path. Moreover, since $Q'$ contains $v_1$ and $y$, and since $h(v_1) = h(y)$, it follows that $Q'$ is a $uv$-NSP.

Now we need to show that if $h(t) = h(s) - k$, then the algorithm above always returns a path. We consider the iteration 
of the algorithm in which $x, y \in V(P_u)$, and $t=v_1$, and $v_1, \dots, v_{k+2} \in V(P_v)$. We claim that the 
subpath $P'$ of $P$ from $v_1$ to $y$ is contained in $H$. Since every vertex $z$ in $V(Q_u) \setminus \sset{x}$ satisfies 
$h(z) \leq h(t) - 2$, it follows from (1) that $z$ has no neighbours in $P'$. Similarly, no vertex in $V(Q_v)$ has a neighbour in $P'$. 
Since $x, y \in V(P_u)$, it follows that the only neighbour of $x$ in $P'$ is $y$. Since $v_1, \dots, v_{k+2} \in V(P_v)$, it follows 
that the only possible neighbour of $v_2, \dots, v_{k+2}$ in $P'$ is $v_1$. This proves our claim. Since $P'$ is a path from $v_1$ to 
$y$ in $H$, it follows that the algorithm returns a path $Q'$. This proves (2).

\bigskip

By (2), we may assume that $h(s)-h(t)\ge 6$.
Let 
$s_0, s_1, \dots, s_6, t_1, \dots, t_6, t_7 \in V(G)$ be distinct, such that: 
\begin{itemize}
\item $s_0\dd s_1\dd s_2\dd s_3$, $s_4\dd s_5\dd s_6$, $t_1\dd t_2\dd t_3$, and $t_4\dd t_5\dd t_6\dd t_7$ are paths;
\item $h(s_i)=h(t_i)$ for $1\le i\le 6$;
\item $h(s_0)+3=h(t_1)+2=h(t_2)+1=h(t_3) \le h(t_4)=h(t_5)-1=h(t_6)-2=h(t_7)-3$;
\item $s_i$ is non-adjacent to $t_j$ for all $i\in \{0\ll 6\}$ and $j\in \{1\ll 7\}$.
\end{itemize}
For each such 14-tuple $s_0, s_1, \dots, s_6, t_1, \dots, t_6, t_7$, we do the following:
\begin{itemize}
    \item We pick a monotone path $Q_u$ from $s_0$ to $u$, and a monotone path $Q_v$ from $t_7$ to $v$. 
    \item We check using \ref{lem:manyshortpaths} whether there are monotone paths $R_u$, $R_v$ such that $R_u$ is an $s_3s_4$-path, 
$R_v$ is a $t_3t_4$-path, and there are no edges between $R_u$ and $R_v$; if not, we move on to the next 14-tuple. 
    \item Let $P_u'$ and $P_v'$ be respectively the concatenations 
$$u\dd Q_u\dd s_0\dd s_1\dd s_2\dd s_3\dd R_u\dd s_4\dd s_5\dd s_6$$
$$t_1\dd t_2\dd t_3\dd R_v\dd t_4\dd t_5\dd t_6\dd t_7\dd Q_v\dd v.$$
Let $H$ be obtained from $G$ by deleting all vertices of $P_u'\setminus \{s_6\}$ and all their neighbours except $s_6$, and deleting
all vertices of $P_v'\setminus \{t_1\}$ and all their neighbours except $t_1$.
We check if there is an induced path $Q$ from $t_1$ to $s_6$ in $H$,
    and if so, we return the concatenated path $u\dd P_u'\dd s_6\dd Q\dd t_1\dd P_v'\dd v$. 
\end{itemize}

If this returns a path $Q'$, then the construction implies that $Q'$ is an induced path; and since $Q'$ contains $s_1, t_1$ with 
$h(s_1)=h(t_1)$, it follows that $Q'$ is a $uv$-NSP. It remains to show that if a shortest $uv$-NSP $P$ exists with 
$h(s) - h(t) \geq 6$, then this algorithm returns a path. We consider the 14-tuple such that 
$s_6 = s$, and $t_1 = t$, $\sset{s_0, \dots, s_6} \subseteq V(P_u)$, and $\sset{t_1, \dots, t_7} \subseteq V(P_v)$. 
This 14-tuple exists since $h(s)-h(t)\ge 6$,
and so there are at least six vertices in $P_u$ that each have the same height as some vertex in $P_v$. 

Now we need to show that the last bullet above returns a path. Let $P'$ be the subpath of $P$ from $s$ to $t$. It follows from (1) that there 
are no edges from $V(Q_u)$ or $V(Q_v)$ to $V(P')$. Since $\sset{s_0, \dots, s_6} \subseteq V(P_u)$ and $\sset{t_1, \dots, t_7} \subseteq V(P_v)$, 
it follows that the only edges from $\sset{s_0, \dots, s_6, t_1, \dots, t_7}$ to $V(P')$ are the edge from $s = s_6$ to its neighbour 
in $V(P')$, and the edge from $t = t_1$ to its neighbour in $V(P')$. If neither $V(R_u)$ nor $V(R_v)$ intersects or has edges to $V(P')$, 
then  $P'$ is present in $H$, and a path is returned. By symmetry, we may assume (for a contradiction)
that $V(R_u)$ intersects or has edges to $V(P')$. 
Let $z$ be the 
vertex closest to $s_3$ in $R_u$ such that $z$ has a neighbour in $V(P')$. 

Let $x\in V(P')$ be the neighbour of $z$ closest to $t=t_1$ in $P'$.
Let $R$ be the induced $uv$-path that begins with a subpath of $P_u'$ from 
$u$ to $z$ and the edge $zx$, and whose remaining vertices are contained in the vertex set of the subpath of $P'$ from $x$ to $t$,
and $P_v'$. Then $R$ is shorter than $P$, since the subpath of $R$ from $u$ to $x$ has length $h(z) + 1$, 
but in $P$, the subpath from $u$ to $x$ contains $s$, and    
thus it has length at least
$h(s) + 1 > h(z) + 1$.
Since $R$ is induced, it follows that 
$R$ is monotone, and therefore $h(x) > h(z)$ (and $x$ has a neighbour in $V(P_v')$, but we will not need this).

The concatenation $Q''$ of the subpath of $P_u'$ from $u$ to $z$, 
the edge $zx$, and the subpath of $P$ from $x$ to $v$ is not monotone, since it contains $s_1$ and $t_1$; and 
as before, it is shorter than $P$.
Therefore $Q''$ is not an induced path. This implies that some vertex $y$ of $P_v$ has a neighbour in
the subpath of $R_u$
between $s_3$ and $z$;  choose $y$ with $h(y)$ maximum, and let $z'$ be a neighbour of $y$ in the subpath of $R_u$
between $s_3$ and $z$, chosen with $h(z')$ maximum (possibly $z'=z$). It follows that $y$ lies in the subpath of $P_v$ between $t_3,t_4$.

Let $t'$ be a vertex of the subpath of $P'$ between $x$ and $t$, such that $h(t')=h(t)$, and subject to that, the subpath of $P'$
between $x,t'$ is minimal.
Now let $R'$ be the concatenation of  a monotone path from $u$ to $t'$, the subpath of $P'$ from $t'$ to $x$, the edge $xz$,
the subpath of $R_u$ between $z$ and $z'$, the edge $z'y$, and the subpath of $P_v$ from $y$ to $v$.
Then $R'$ is an induced path because of (1); and its 
length is at most the length of $P'$ plus $d(u,t) + 2 + d(y,v)$; but the length of $P$ is at least the length of $P'$ plus 
$d(u, t) + 6 + d(t, v)$, 
and $d(t, v) \geq d(y,v)$ since $y \in V(P_v)$. This implies that $R'$ is monotone. Since $z$ is closer to $v$ than $x$ in $R'$, 
it follows that $h(x) < h(z)$, a contradiction. Hence the last bullet above does indeed return a path. 
(We omit the analysis of running time, which is straightforward.) This proves \ref{thm:k0}.~\bbox

Now \ref{thm:long} follows from \ref{lem:components} and \ref{thm:k0}.

\section{Acknowledgments}

The first author was supported by Israel Science Foundation Grant 100004639 and Binational Science Foundation USA-Israel Grant 100005728. The second author was supported by AFOSR grant A9550-19-1-0187 and NSF
grant DMS-1800053. This material is based upon work supported by the National Science Foundation under Award No. DMS-1802201 (Spirkl).

\end{document}